\documentclass{amsart}

\usepackage{hyperref}

\newtheorem{thm}{Theorem}[section]

\newtheorem{prop}[thm]{Proposition}
\newtheorem{lem}[thm]{Lemma}
\theoremstyle{definition}

\numberwithin{equation}{section}

\newcommand{\Hom}{\mbox{Hom}}
\newcommand{\Spec}{\mbox{Spec}}

\begin{document}  

\subjclass{Primary 14H45. Secondary 14E22, 14M25.}
\title{The Genus of a Curve of Fermat Type} 
\author{Jeremiah M. Kermes} 
\address{5022 Holly Ridge Dr.\\ Raleigh, NC 27612}
\email{cermz27@bellsouth.net} 

\begin{abstract}
In this paper we begin to study curves on a weighted projective plane with one trivial weight, ${\mathbb P}(1,m,n)$, by determining the genus of curves of Fermat type.  These are curves, $C$, defined by the ``homogeneous'' polynomial $x_0^{amn}+x_1^{an}=x_2^{am}$.  We begin by finding local coordinates for the standard affine cover of ${\mathbb P}(1,m,n)$, and then prove that the curve is smooth.  This is done by pulling the curve up to the surface's desingularization, ${\mathbb D}(1,m,n)$.  Then a map, $\phi:C\rightarrow{\mathbb P}^1$ is constructed, and it's ramification divisor is determined.  We conclude by applying Hurwitz's theorem to $\phi$ to obtain $C$'s genus.  
\end{abstract}

\maketitle

We begin by making some simplifications.  First, a result of Dolgachev \cite[1.3.1]{Dol} tells us that ${\mathbb P}(1,am,an)\cong {\mathbb P}(1,m,n)$.  Thus, we may assume that $m$ and $n$ are relatively prime.  Next, note that $[x_0,x_1,x_2]\mapsto [x_0,x_2,x_1]$ gives an isomorphism of the coordinate rings for ${\mathbb P}(1,m,n)$ and ${\mathbb P}(1,n,m)$.  Subsequenly, we may assume that $m<n$.  

In section \ref{sec:classical} an affine cover for ${\mathbb P}(1,m,n)$ is constructed using the classical description of a weighted projective space as $Proj$ of a graded ring.  The local (affine) forms of the equations defining a Fermat-type curve are then constructed.  Section \ref{sec:toric} gives the construction of ${\mathbb P}(1,m,n)$ as a toric variety and describes the isomorphism with the classical construction by relating the generators of the rings for the affine cover.  The goal of section \ref{sec:ideals} is to use the machinery of toric varieties write to down some of the polynomials in the monomial ideal ${\mathcal I}(U_1)$ (resp. ${\mathcal I}(U_2)$) of functions vanishing on the affine surface $U_1$ (resp. $U_2$).  In section \ref{sec:smooth} the desingularization algorithm in \cite{Ker} and \cite{Oda} is used to show that a Fermat-type curve is smooth, while section \ref{sec:ramification} constructs a map from the curve to ${\mathbb P}^1$ and the ramifacation divisor of this map.  Finally, Hurwitz's theorem is used to determine the genus of a Fermat-type curve.  Throughout this paper $k$ is an algebraically closed field of characteristic $0$.

\section{The Classical Description} 
\label{sec:classical}
The standard affine cover of a weighted projective plane consists of three affine varieties, $U_i=\{[x_0,x_1,x_2]\in{\mathbb P}(a_0,a_1,a_2)|$ $x_i\neq0\}$.  This is the space $Spec(k[U_i])$ where $k[U_i]$ is the degree zero part of the graded ring $k[x_0,x_1,x_2][x_i^{-1}]$ where the grading is given by $deg(x_i)=a_i$.  

In the case of ${\mathbb P}(1,m,n)$ we can use the fact that $deg(x_0)=1$ to see that $U_0\cong{\mathbb A}^2$.  This is because the degree zero part of $k[x_0,x_1,x_2][x_0^{-1}]$ is generated by $\frac{x_1}{x_0^m}$ and $\frac{x_2}{x_0^n}$.  Since there are no relations on these forms we see $U_0=Spec\left(k\left[\frac{x_1}{x_0^m},\frac{x_2}{x_0^n}\right]\right)=Spec\left(k\left[X,Y\right]\right)={\mathbb A}^2$.  

The other two affine surfaces, however, are singular.  Begin by finding generators for the $k$-algebras, $k[U_1]$ and $k[U_2]$.  On $U_1$, the polynomial generator $x_1$ becomes a unit.  Thus, $k[x_0,x_1,x_2][x_1^{-1}]$ is generated over its quotient field by $x_0$ and $x_2$.  

\begin{lem}:$\:$\\
\label{lem:U1gen}
$k[U_1]$ is generated as an affine $k$-algebra by \begin{equation}\label{eq:U1gen} z_j=\frac{x_0^{m\left\lceil\frac{jn}{m}\right\rceil-jn}x_2^j}{x_1^{\left\lceil\frac{jn}{m}\right\rceil}}\end{equation} where $0\leq j\leq m$ and $\left\lceil x\right\rceil$ is the next largest integer than $x$, with the exception that $\left\lceil0\right\rceil=1$.  
\end{lem}
\begin{proof}
The goal is to generate the forms of degree $0$ where $x_1$ is invertible.  Begin with the obvious form, $z_0=\frac{x_0^m}{x_1}$, which is the reason for taking $\left\lceil0\right\rceil=1$ to avoid a useless generator of $1$.  

The technique to generate the remaining $z_j$'s is to take successively higher powers of $x_2$.  Then put just enough $x_1$ terms in the denominator to make the degree of the whole form negative, which is to say $\left\lceil\frac{jn}{m}\right\rceil$ of them.  Then to give a form of degree zero, take advantage of the fact that $deg(x_0)=1$ and put just enough $x_0$ terms in the numerator to give the form a total degree of zero.  The number of $x_0$'s necessary is then \[ \left\lceil\frac{jn}{m}\right\rceil\cdot deg(x_1)-j\cdot deg(x_2)=m\left\lceil\frac{jn}{m}\right\rceil-jn .\]  With this method one keeps generating possibly distinct forms until the final form of $z_m=\frac{x_2^m}{x_1^n}$ is reached.  
\end{proof}

In fact the exact same technique with the roles of $x_1$ and $x_2$ reversed will yield the corresponding result for the other singular affine surface, $U_2$.  By continuing with the convention that $\left\lceil0\right\rceil=1$ we have the following Lemma.  

\begin{lem}:$\:$\\
\label{lem:U2gen}
$k[U_2]$ is generated as an affine $k$-algebra by \begin{equation}\label{eq:U2gen} w_j=\frac{x_0^{n\left\lceil\frac{jm}{n}\right\rceil-jm}x_1^j}{x_2^{\left\lceil\frac{jm}{n}\right\rceil}} \end{equation} where $0\leq j\leq n$.  
\end{lem}

Now that we have the rings for an affine cover of ${\mathbb P}(1,m,n)$ we can write down the equation for a Fermat-type curve locally. Recall that such a curve is written in terms of it homogeneous coordinate ring as \begin{equation}\label{eq:ferdef} x_0^{amn}+x_1^{an}-x_2^{am}.\end{equation}  On $U_0$ the element $x_0$ is a unit so that this equation becomes $1+X^{an}-Y^{am}$.  By a similar process on $U_1$ using the local coordinates $\left(z_0,\ldots,z_m\right)$ we see that the curve is defined by \begin{equation}\label{eq:WFC-U1} 0=\left(\frac{x_0^m}{x_1}\right)^{an}+1-\left(\frac{x_2^m}{x_1^n}\right)^a=z_0^{an}+1-z_m^a.\end{equation}  The corresponding result for $U_2$ is then $w_0^{am}+w_n^a-1$.  

\section{The Toric Approach}
\label{sec:toric}

In this paper we use the description and notation for toric varieties found in \cite{Fult}

A weighted projective space \cite{Dol}, ${\mathbb P}(a_0,\ldots,a_d)$, is the complete toric variety whose fan is given by $\Delta(1)=\left\{v_0,\ldots, v_d\right\}$ where $Span_{\mathbb Z}(v_0,\ldots,v_d)={\mathbb Z}^d={\mathbb N}$ is the lattice, and $\sum_{\Delta(1)} a_jv_j=0$.  In the case of ${\mathbb P}(1,m,n)$ such a fan can be given by \[ \Delta(1)=\left\{\begin{bmatrix}-m\\ -n\end{bmatrix}, \begin{bmatrix}1\\ 0\end{bmatrix}, \begin{bmatrix}0\\ 1\end{bmatrix}\right\} \]  where each maximal cone $\sigma_i$ is the cone generated over ${\mathbb R}_+$ by $\Delta(1)\setminus\{v_i\}$.  

A quick check of the determinants of the edges of each maximal cone will show that the surfaces $U_{\sigma_1}$ and $U_{\sigma_2}$ will be singular, while $U_{\sigma_0}$ is smooth.  To describe $k[U_0]$ note that $\sigma_0=\left<e_1,e_2\right>$ where $e_1,e_2$ is the standard basis for ${\mathbb Z}^2$.  Then the dual cone is $\sigma_0^\vee=\left<e_2^\vee,e_1^\vee\right>$ where $e_1^\vee, e_2^\vee$ is dual to the standard basis.  Then taking $X=\chi^{e_1^\vee}$ and $Y=\chi^{e_2^\vee}$ gives $k[U_{\sigma_0}]\cong k[Y,X]$.  

Since the remaining surfaces are singular, it will be more difficult to construct their affine rings.  To begin, note that their dual cones are \[ \sigma_1^\vee=\left<\begin{bmatrix}-1\\ 0\end{bmatrix}, \begin{bmatrix}-n\\ m\end{bmatrix}\right>\quad \mbox{ and } \quad \sigma_2^\vee=\left<\begin{bmatrix}n\\ -m\end{bmatrix}, \begin{bmatrix}0\\ -1\end{bmatrix}\right>\] respectively.  

Proceed with finding generators of the semi-group $S_{\sigma_1}=\sigma_1^\vee\cap{\mathbb N}^*$ by finding lattice points in the parallelogram formed by $-e_1^\vee$ and $-ne_1^\vee+me_2^\vee$.  Do this by considering rational numbers $0\leq s,t\leq1$ where $s=0$ if and only if $t=1$ and vice-versa.  Then the generators of $S_{\sigma_1}$ are lattice points of the form $s\left(-e_1^\vee\right)+t\left(-ne_1^\vee+me_2^\vee\right)=-\left(s+tn\right)e_1^\vee+\left(tm\right)e_2^\vee$.  For the second coordinate to be integral we must have $t=\frac{j}{m}$ for $0\leq j\leq m$.  Since $0\leq s\leq1$ the first coefficient must then be $\left\lceil\frac{nj}{m}\right\rceil$ where $\left\lceil\cdot\right\rceil$ denotes the next largest integer.  It should be noted that in order to match the condition $t=0\Rightarrow s=1$ we must again use the convention that $\left\lceil0\right\rceil=1$.  

Adopting this convention we see that generators for $S_{\sigma_1}$ are the lattice points $u_j=-\left\lceil\frac{nj}{m}\right\rceil e_1^\vee+je_2^\vee$ for $0\leq j\leq m$.  The generators for the $k$-algebra, $k[\sigma_1]$, are then \begin{equation}\label{eq:U1torcoord} \tilde{z}_j=\chi^{u_j}=X^{-\left\lceil\frac{nj}{m}\right\rceil}Y^j. \end{equation}  By a similar argument one sees that $k[\sigma_2]$ is generated by \begin{equation}\label{eq:U2torcoord} \tilde{w}_j=X^jY^{-\left\lceil\frac{mj}{n}\right\rceil}\end{equation} for $0\leq j\leq n$.  

\begin{thm}\label{thm:tor-class-iso}
The isomorphism between classical and toric constructions of ${\mathbb P}(1,m,n)$ is given by $X\sim\frac{x_1}{x_0^m}$ and $Y\sim\frac{x_2}{x_0^n}$ where $X=\chi^{e_1^\vee}$ and $Y=\chi^{e_2^\vee}$.  
\end{thm}
\begin{proof}
We prove this by showing that this correspondence gives an isomorphism of the three $k$-algebras giving the affine cover of the surface in each construction.  In particular, it gives an isomorphism on the generators of these rings such that $U_i\cong U_{\sigma_i}$.  

Checking this on $U_0$ is trivial since $k\left[Y,X\right]\cong k\left[\frac{x_2}{x_0^n},\frac{x_1}{x_0^m}\right]$.  To see the isomorphism on $U_1$ note that the correspondence means \[ \tilde{z}_j=X^{-\left\lceil\frac{nj}{m}\right\rceil}Y^j=\left(\frac{x_1}{x_0^m}\right)^{-\left\lceil\frac{nj}{m}\right\rceil}\left(\frac{x_2}{x_0}\right)^j=\frac{x_0^{m\left\lceil\frac{nj}{m}\right\rceil-nj}x_2^j}{x_1^{\left\lceil\frac{nj}{m}\right\rceil}}=z_j\] while a similar calculation gives $\tilde{w}_j=w_j$ to complete the isomorphism on $U_2$.  
\end{proof}

\section{The Ideals for the Affine Cover}
\label{sec:ideals}

The homomorphism $k[z_0,\ldots,z_m]\rightarrow k[\sigma_1]$ gives an embedding $U_{\sigma_1}\hookrightarrow{\mathbb A}^{m+1}$.  The kernel of this homomorphism, ${\mathcal I}(U_1)$, is the ideal of functions on ${\mathbb A}^{m+1}$ vanishing on $U_1$.  In order to determine some of the equations in ${\mathcal I}(U_1)$ note that $k[\sigma_1]$ is generated by $z_j=\chi^{u_j}$ for $0\leq j\leq m$ where the $u_j$'s are lattice points in the cone $\sigma_1^\vee$.  By regarding a linearly independent pair $u_{i_1}$, $u_{i_2}$  as a basis for the vector space ${\mathbb N}^*\otimes{\mathbb Q}$, one can write the remaining $u_j\in S_{\sigma_1}$ as a rational linear combination of them.  Multiplication by the common denominator of these fractions yields an equation with integral coefficients $a_ju_j=b_ju_{i_1}+c_ju_{i_2}$ that the map $u_j\mapsto\chi^{u_j}$ turns into \[ z_j^{a_j}=z_{i_1}^{b_j}z_{i_2}^{c_j}.\]  Any negative exponents can be multiplied out to convert this to a polynomial in ${\mathcal I}(U_1)$.  For convenience we refer to a set of polynomials obtained in this fashion as having type $(i_1,i_2)$.  

For example, the type $(0,m)$ equations are nearly already done for us since $u_0$ and $u_m$ are the edges of $\sigma_1^\vee$ that were used to find the other generators of $S_{\sigma_1}$.  Recall from the derivation of equation \ref{eq:U1torcoord} that $u_j=su_0+tu_m$ where $t=\frac{j}{m}$ and $s=\left\lceil\frac{nj}{m}\right\rceil-\frac{nj}{m}$.  Multiplying both sides by $m$ and applying $\chi$ yields the $(0,m)$ polynomials of ${\mathcal I}(U_1)$ as \begin{equation}\label{eq:U1-0,m} z_j^m-z_0^{m\left\lceil\frac{nj}{m}\right\rceil-nj}z_m^j \end{equation} for $1\leq j\leq m-1$  Using the same technique shows the $(0,n)$ equations of ${\mathcal I}(U_2)$ to be \begin{equation}\label{eq:U2-0,n} w_j^n=w_0^{n\left\lceil\frac{mj}{n}\right\rceil-mj}w_n^j \end{equation} for $1\leq j\leq n-1$.  

It is not difficult to use this same technique to generate the type $(0,1)$ equations as well.  Simply note that $u_0=-e_1^\vee$ and $u_1=-\left\lceil\frac{n}{m}\right\rceil e_1^\vee+e_2^\vee$.  The result in this case is that ${\mathcal I}(U_1)$ will contain \begin{equation}\label{eq:U1-0,1} z_1^j-z_0^{j\left\lceil\frac{n}{m}\right\rceil-\left\lceil\frac{nj}{m}\right\rceil}z_j\end{equation} for $2\leq j\leq m$.  Similiarly, since $0<m<n$, ${\mathcal I}(U_2)$ contains \begin{equation}\label{eq:U2-0,1} w_1^j-w_0^{j-\left\lceil\frac{mj}{n}\right\rceil}w_j\end{equation} for $2\leq j\leq n$.  

\section{An Aside on the Smoothness of Fermat-type Curves}
\label{sec:smooth}

Showing that a Fermat-type curve is smooth will allow certain machinery to be applied to the study of these curves.  Begin by noting that on $U_0$, the curve $C$ is defined by \[ 1+\left(\frac{x_1}{x_0^m}\right)^{an}-\left(\frac{x_2}{x_0^n}\right)^{am}.\]  In terms of the toric coordinates, $X=\chi^{e_1^\vee}$ and $Y=\chi^{e_2^\vee}$, this is $1+X^{an}-Y^{am}$.  Plugging this into the Jacobi criterion will show that $C$ is smooth where $x_0\neq0$.  

All that remains now is to check that $C$ is smooth where it intersects the line $x_0=0$.  Since the curve must satisfy $x_0^{amn}+x_1^{an}-x_2^{am}=0$ any points on this line will satisfy $x_1^{an}=x_2^{am}$.  Consequently $x_1=0 \Leftrightarrow x_2=0$.  This means that any points of $C$ with $x_0=0$ are contained in $U_1\cap U_2$.  Thus we need only verify the smoothness of $C$ on $U_1$.  

To do this one could try combining Equations \ref{eq:U1-0,m} and \ref{eq:U1-0,1} with Equation \ref{eq:WFC-U1} and using the Jacobi criterion.  This attempt, however, would fail to do anything more than show that possible singular points lie on the line $x_0=0$, which we already know.  This is because Equations \ref{eq:U1-0,m} and \ref{eq:U1-0,1} fail to generate the whole ideal ${\mathcal I}(U_1)$.  

On the other hand, because ${\mathbb P}(1,m,n)$ is a complete toric variety, it is normal so that any singularities have codimension at least $2$.  Since this is a surface the singularities will be isolated to the fixed points of the toric action on $U_1$ and $U_2$, which are $[0,0,1]$ and $[0,1,0]$ in terms of homogeneous coordinates.  But $C$ contains neither of these points.  Thus, if we consider this surface's desingularization $\pi:{\mathbb D}(1,m,n)\rightarrow{\mathbb P}(1,m,n)$ as in \cite{Ker}, we obtain an isomorphism $C\cong\pi^{-1}(C)$.  

In fact the concern here is not the entirety of ${\mathbb D}(1,m,n)$, but rather $\pi^{-1}(U_1)$ since we merely have to demonstrate the smoothness of $C$ on $U_1$.  It is known that the desingularization of an affine toric surface corresponding to a cone, $\sigma$ is the toric surface obtained by subdividing $\sigma$ throught the rays $\left\{l_0,\ldots,l_{s+1}\right\}$ given by \cite[Lemma 1.20]{Oda} where $l_0$ and $l_{s+1}$ are the edges of $\sigma$.  The maximal cones for $\pi^{-1}(U_\sigma)$ are then $\tau_j=\left<l_{j-1},l_j\right>$.

An additional part of the algorithm of great importance is a collection of integers $\left\{b_1,\ldots,b_s\right\}$ with each $b_j\geq2$.  Geometrically, these numbers correspond to the self-intersection number of the $T$-equivariant divisors on $\pi^{-1}(U_\sigma)$ by $D(l_j)=-b_j$.  It is also important to note from \cite[Prop. 1.19]{Oda} that they satisfy \begin{equation}\label{eq:Oda} l_{j-1}+l_{j+1}=b_jl_j. \end{equation}

We define elements of the dual lattice ${\mathbb M}=\Hom({\mathbb N},{\mathbb Z})$ by letting $l_j^\perp$ be the unique element with $l_j^\perp(l_j)=0$ and $l_j^\perp(l_{j-1})=1$ (or equivalently $l_j^\perp(l_{j+1})=-1$).  This is well defined since $\tau_j$ is non-singular, so $\det|l_{j-1},l_j|=\pm1$.  

Using this notation each of the $k$-algebras $k[\tau_j]$ is simply $k[x_j,y_j]$ where $x_j=\chi^{-l_{j-1}^\perp}$ and $y_j=\chi^{l_j^\perp}$.  Next up, we need to know how to change coordinates between $U_{\tau_j}$ and $U_{\tau_{j+1}}$.  

\begin{lem}\label{lem:coord-change}
The $k$-algebra isomorphism between $k[\tau_j][y_j^{-1}]$ and $k[\tau_{j+1}][x_{j+1}^{-1}]$ is given by $x_j\mapsto x_{j+1}^{b_j}y_{j+1}$ and $y_j\mapsto x_{j+1}^{-1}$.
\end{lem}
\begin{proof}
Both of these algebras are simply $k[\tau_j\cap\tau_{j+1}]$.  In this region we may invert the element corresponding to their common edge, $l_j$.  This leads one to observe that $y_j=\chi^{l_j^\perp}=\left(\chi^{-l_j^\perp}\right)^{-1}=x_{j+1}^{-1}$.  

To prove the $x_j$ piece of the isomorphism note that it is equivalent to the statement $-l_{j-1}^\perp=-b_jl_j^\perp+l_{j+1}^\perp$.  This will be proven by showing that $l_{j-1}^\perp+l_{j+1}^\perp-b_jl_j^\perp$ vanishes on a basis for ${\mathbb N}={\mathbb Z}^2$ (and hence, on all of ${\mathbb N}$).  Since $U_{\tau_j}$ is smooth the vectors $l_{j-1}$ and $l_j$ constitute a suitable basis.  

Begin by recalling that $l_j^\perp(l_j)=0$, $l_{j+1}^\perp(l_j)=1$, and $l_{j-1}^\perp(l_j)=-1$.  Subsequently we have $l_{j-1}^\perp(l_j)+l_{j+1}^\perp(l_j)-b_jl_j^\perp(l_j)=-1+1-0$ to show that the form vanishes on $l_j$.  Using the same process for $l_{j-1}$ gives $l_{j+1}^\perp(l_{j-1})-b_j$.  In order to determine $l_{j+1}^\perp(l_{j-1})$ solve equation \ref{eq:Oda} for $l_{j-1}$ and use the linearity of $l_{j+1}^\perp$ to obtain \[l_{j+1}^\perp(l_{j-1})=b_jl_{j+1}^\perp(l_j)-l_{j+1}^\perp(l_{j+1})=b_j \]  which can be plugged back in to find $l_{j-1}^\perp(l_{j-1})+l_{j+1}^\perp(l_{j-1})-b_jl_j^\perp(l_{j-1})=b_j-b_j=0$, concluding the proof.  
\end{proof}

This isomorphism will allow us to write the polynomial defining a Fermat-type curve on each of the open affine neighborhoods $U_{\tau_j}=\mbox{Spec}\left(k[x_j,y_j]\right)$.  In order to complete this process, a couple of auxiliary sequences will need to be obtained.  

The first of these sequences, $\left\{r_{-1},\ldots,r_{s-1}\right\}$, was constructed in \cite[Theorem 6.1]{Ker}.  In the case of $\sigma_1$ for ${\mathbb P}(1,m,n)$, by expressing $n$ as $mk+r$ with $m,r$ relatively prime, the inital values in this sequence are $r_{-1}=m$ and $r_0=r$.  It was also shown that this is a sequence of positive integers satisfying $r_j=b_jr_{j-1}-r_{j-2}$.  

The other sequence, $\left\{t_0,\ldots,t_s\right\}$  is given by $t_0=0$, $t_1=1$ and $t_{j+1}=b_jt_j-t_{j-1}$.  Using the fact that every $b_j\geq2$, it is a simple matter to prove inductively that this sequence is increasing, which means each $t_j$ is non-negative.  

\begin{lem}
\label{lem:localWFC} 
On the region $U_{\tau_j}$ a Fermat-type curve of degree $amn$ is determined by the polynomial \begin{equation}\label{eq:localWFC} F_j=x_j^{an\cdot t_j}y_j^{an\cdot t_{j-1}}-x_j^{a\cdot r_{j-1}}y_j^{a\cdot r_{j-2}}+1\end{equation} for $1\leq j\leq s+1$.  
\end{lem}
\begin{proof}
The proof is by induction on $j$.  When $j=1$, \cite[Lemma 4.1]{Ker} shows that $\tau_1=\left<\begin{bmatrix}0\\1\end{bmatrix},\begin{bmatrix}-1\\ -k\end{bmatrix}\right>$.  This means that in terms of $X=\chi^{e_1^\vee}$ and $Y=\chi^{e_2^\vee}$ the local coordinates for $U_{\tau_1}$ are $x_1=X^{-1}$ and $y_1=X^{-k}Y$.  Recall from equation \ref{eq:WFC-U1} that the curve on $U_{\sigma_1}$ is defined by $z_0^{an}-z_m^a+1$.  Theorem \ref{thm:tor-class-iso} allows us to write this in toric coordinates as $X^{-an}-X^{-an}Y^{am}+1$.  A little algebra gives the corresponding polynomial on $U_{\tau_1}$ to be $x_1^{an}-x_1^{ar}y_1^{am}+1$ proving the case $j=1$.  

Now suppose the Lemma holds for $F_j$.  Then $F_{j+1}$ will be the image of $F_j$ under the map from Lemma \ref{lem:coord-change}.  This is turns out to be \[ \left(x_{j+1}^{b_j}y_{j+1}\right)^{an\cdot t_j}\left(x_{j+1}^{-1}\right)^{an\cdot t_{j-1}}-\left(x_{j+1}^{b_j}y_{j+1}\right)^{a\cdot r_{j-1}}\left(x_{j+1}^{-1}\right)^{a\cdot r_{j-2}}+1. \]  Collecting terms and using the recursive definition of $r_j$ and $t_{j+1}$ then gives \[ x_{j+1}^{an\cdot t_{j+1}}y_{j+1}^{an\cdot t_j}-x_{j+1}^{a\cdot r_j}y_{j+1}^{a\cdot r_{j-1}}+1 \] concluding the proof.  
\end{proof}

With the various local formulations of $\pi^{-1}(C\cap U_1)$ in hand it is not difficult to prove the following.  

\begin{thm}
\label{thm:WFCsmooth}
A Fermat-type curve is smooth.
\end{thm}
\begin{proof}
The only part of the proof that remains unfinished is checking the points of $C$ with $x_0=0$, which all lie in $U_1$.  Since $C\cong\pi^{-1}(C)$ this means we just need to check for smoothness on $\pi^{-1}(U_1)=\bigcup_{j=1}^{s+1}U_{\tau_j}$.  

Begin by noting that for the interior cones $2\leq j\leq s$ any point of $\pi^{-1}(C)\cap U_{\tau_j}$ satisfies $x_j\neq0$ and $y_j\neq0$.  Thus each such point is also contained in $U_{\tau_{j-1}}$, so we really only need to check the cases $j=1$ and $j=s+1$.  

On $U_{\tau_1}$ note that $C$ is defined by $x_1^{an}-x_1^{ar}y_1^{am}+1$ which contains no points with $x_1=0$.  Also, the differential is \[dF_1=\left(an\cdot x_1^{an-1}-ar\cdot x_1^{ar-1}y_1^{am}\right)dx_1-\left(am\cdot x_1^{ar}y_1^{am-1}\right)dy_1.\]  Since $x_1\neq0$ the only way for the $dy_1$ coefficient to vanish is to have $y_1=0$.  However, this results in a non-zero $dx_1$ coefficient, so that the curve is smooth on $U_{\tau_1}$.  

To handle the case $j=s+1$ one must recall the sequence of rational numbers $\left\{\beta_0,\ldots,\beta_{s-1}\right\}$ in \cite[Eq. 2.1]{Ker} which are related to the $r_j$'s by $\beta_j=\frac{r_{j-1}}{r_j}$.  In particular, \cite[Lem 2.1]{Ker} tells us that $\frac{r_j}{r_{j-1}}=b_j-\beta_{j-1}$.  Since the last $\beta_j$ occurs when $\beta_{s-1}=b_s\in{\mathbb Z}$ this means that $\frac{r_s}{r_{s-1}}=b_s-\beta_{s-1}=0$.  Consequently, $r_s=0$ (and $r_{s-1}>0$).  Since the $t_j$'s are increasing this leaves the polynomial \[ F_{s+1}=x_{s+1}^{an\cdot t_{s+1}}y_{s+1}^{an\cdot t_s}-x_{s+1}^{a\cdot r_{s-1}}+1 \] whose resulting curve contains no points with $x_{s+1}=0$.  Applying the Jacobi criterion as was done in the $j=1$ case will complete the proof.
\end{proof}

\section{Ramification of a Map to the Projective Line}
\label{sec:ramification}

Now that we know a Fermat-type curve is smooth we can proceed with determining its genus.  The approach will be to construct a map to ${\mathbb P}^1$, determine the degree and ramification divisor of this map, and at last use Hurwitz's Theorem \cite[Cor. IV.2.4]{Hart} to determine $C$'s genus.  

The map used will be the rational map $\phi:{\mathbb P}(1,m,n)\rightarrow{\mathbb P}(1,m)$ sending $[x_0,x_1,x_2]$ to $[x_0,x_1]$ in terms of homogeneous coordinates.  The only point at which this is undefined, $[0,0,1]$ is not on the curve, so it restricts to a morhpism on $C$.  The following Lemma reveals that ${\mathbb P}(1,m)$ is a rather simple space.  

\begin{prop}\label{prop:WPline}
If $s,t$ are relatively prime, then ${\mathbb P}(s,t)\cong{\mathbb P}^1$.  
\end{prop}
\begin{proof}
Since ${\mathbb P}(s,t)=\mbox{Proj}\left(k[x_0,x_1]\right)$ with the grading $\deg(x_0)=s$, $\deg(x_1)=t$, it is covered by the two regions $V_i=\left\{[x_0,x_1]\mbox{ s.t. }x_i\neq0\right\}$.  Specifically, \begin{equation}\label{eq:P1local} V_0=\Spec\left(k\left[\frac{x_1^s}{x_0^t}\right]\right)\quad V_1=\Spec\left(k\left[\frac{x_0^t}{x_1^s}\right]\right).\end{equation} The resulting space is two affine lines with a coordinate change $x\mapsto x^{-1}$, i.e. ${\mathbb P}^1$.  
\end{proof}

With this map in hand, and the local information obtained in section \ref{sec:ideals}, we can determine both the the degree of $\phi$ and its ramification divisor, $R$.  

\begin{lem}\label{lem:degphi}
Let $C$ be a Fermat-type curve on ${\mathbb P}(1,m,n)$ of degree $amn$.  The degree of the map $\phi:C\rightarrow{\mathbb P}(1,m)\cong{\mathbb P}^1$ given by $[x_0,x_1,x_2]\mapsto [x_0,x_1]$ is $am$.  
\end{lem}
\begin{proof}
The degree will be determined by finding the number of distinct points in a generic fiber of $\phi$.  This may be done on the dense, open subset, $U_0=\Spec\left(k[X,Y]\right)$ where $X=\frac{x_1}{x_0^m}$ and $Y=\frac{x_2}{x_0^n}$.  Using this notation, $k[V_0]=k[X]$, and $\phi|_{U_0}$ corresponds to the inclusion of $k$-algebras, $k[X]\hookrightarrow k[X,Y]$.  Consequently, the restriction of $\phi$ to $C\cap U_0$ is obtained by composing this with the natural projection to \[ k[C\cap U_0]\cong k[X,Y]/\left<1+X^{an}-Y^{am}\right>.\]  Now note that unless $1+X^{an}=0$ (which only happens for finitely many points on $V_0$), this quantity will have $am$ distinct $am^{th}$ roots.  Each of these corresponds to a distinct $Y$-value, yielding an equal number of distinct points in the fiber of $X$ and proving the Lemma.  
\end{proof}

In fact, the points on $V_0$ where $1+X^{an}=0$ are more than just the points where $\phi:C\rightarrow{\mathbb P}^1$ is not $am$-to-one.  These are some of the branch points whose fibers will consists of ramification points of $\phi$.  They are not, as we shall see, all of the branch points.  Thus we begin our determination of the ramification divisor by splitting it into two pieces.  The first, $\bar{R}$, will consist of those ramification points contained in $U_0$.  The other piece, $R_0$ is merely those ramification points on the line $x_0=0$.  

\begin{lem}\label{lem:U0ramdiv}
The ramification divisor for the map $\phi:C\cap U_0\rightarrow{\mathbb P}^1$ is \[ \bar{R}=\sum_{j=1}^{an}\left(am-1\right)\left[1,\alpha_j,0\right]\] where each $\alpha_j$ is a distinct $an^{th}$ root of $-1$.  
\end{lem}
\begin{proof}
Note that since $x_0\neq0$ on both regions, that the image of $U_0$ is wholly contained in $V_0$.  Then the proof of Lemma \ref{lem:degphi} revealed the branch points to be those with $1+X^{an}=0$, which are the $an$ distinct ideals $\left<X-\alpha_j\right>$.  A consequence of this is that the fibers of these branch points must satisfy $Y^{am}=1+\alpha_j^{an}=0$, so the branch points, $X=\alpha_j$, are in one to one correspondence with the ramification points $X=\alpha_j$, $Y=0$ (or in terms of homogeneous coordinates the branch points are $\left[1,\alpha_j\right]$ and the ramifiaction points are$\left[1,\alpha_j,0\right]$).  

Next it remains to find the coefficient of each ramification point $P\in\bar{R}$.  By \cite[Prop. IV.2.2]{Hart} this is merely $v_P(t)-1$ where $t$ is the generator for the one-dimensional maximal ideal $\phi(P)\subset k[X]$ and $v_P$ is the standard valuation at the ramification point.  Now $P$ is the ideal $\left<X-\alpha_j,Y\right>$ and $\phi(P)$ is $\left<t\right>=\left<X-\alpha_j\right>$.  In particular, the valuation, $v_P$, is taking place in the ring $k[C]=k[X,Y]/\left<1+X^{an}-Y^{am}\right>$ localized at the ideal $P=\left<X-\alpha_j,Y\right>$.  Since $C$ is a smooth curve, this point must be a principal ideal generated by either $X-\alpha_j$ or $Y$.  

Since $P\in C$, in ${\mathcal O}_P$ we have $Y^{am}=1+X^{an}=\left(X-\alpha_j\right)\cdot p(X)$ where $p(\alpha_j)\neq0$, so $p(X)$ is a unit.  Subsequently, $X-\alpha_j\in \left<Y\right>$ so that $Y$ generates the principal ideal, $P$.  Furthermore $v_P(X-\alpha_j)=am$ so that the coefficient of $P=[1,\alpha_j,0]$ is $(am-1)$.  
\end{proof}

Now we turn our attention to the ramification points of $C$ that lie on the line $x_0=0$.  We saw in the beginning of section \ref{sec:smooth} that all such points lie in $U_1\cap U_2$, so that we may work in the affine region $U_1$.  This will mean using the local coordinates of Lemma \ref{lem:U1gen}, $\left(z_0,\ldots,z_m\right)$.  Furthermore, these points are all in the fiber of $\phi$ sitting over the origin of $V_1$.  As it turns out this is a branch point of $\phi$ and every point of $C$ with $x_0=0$ is a ramification point.  

\begin{lem}\label{lem:x0ramdiv}
In terms of local coordinates on $U_1$, the ramification occuring on the line $x_0=0$ is \[ R_0=\sum_{j=1}^a\left(m-1\right)\cdot\left(0,\ldots,0,\gamma_j\right)\] where each $\gamma_j$ is a distinct $a^{th}$ rooth of $1$.
\end{lem}
\begin{proof}
Begin by finding all of the points on $C$ with $x_0=0$.  To do this, note that Lemma \ref{lem:U1gen} implies that they are of the form $(0,\ldots,0,z_m)$.  Furthermore equation \ref{eq:WFC-U1} says that any such points on $C$ must satisfy $0=z_0^{an}+1-z_m^a$.  Since $z_0=0$ there are $a$ distinct points on $\phi^{-1}\left([0,1]\right)=C\cap\left\{x_0=0\right\}$ given by the distinct $a^{th}$ roots of unity, $z_m=\gamma_j$.  

Also note that $k[V_1]=k\left[\frac{x_0^m}{x_1}\right]=k\left[X^{-1}\right]=k[z_0]$, and the $k$-algebra homomorphism corresponding to $\phi$ is the map $k[z_0]\rightarrow k[z_0,\ldots,z_m]/{\mathcal I}(U_1)$ sending $z_0$ to itself.  Since all of the points  $P=\left<z_0,\ldots,z_{m-1},z_m-\gamma_j\right>$ project to the origin of $V_1$, $\left<z_0\right>$, all that remains is to show these are ramifiaction points with $v_P\left(z_0\right)=m$.  

To do this note that $z_m(P)=\gamma_j\neq0$ means $z_m\not\in P$ is a unit in the local ring ${\mathcal O}_P$, yielding $v_P(z_m)=0$.  Now consider equation \ref{eq:U1-0,m} with $j=1$, i.e. $z_1^m=z_mz_0^{m-r}$.  Taking valuations of both sides gives \begin{equation}\label{eq:valrel-z0z1} mv_P\left(z_1\right)=\left(m-r\right)v_P\left(z_0\right).\end{equation}  Since $m$ and $m-r$ are relatively prime this means that $m|v_P\left(z_0\right)$.  Now if we can show that $v_P\left(z_0\right)|m$, we'll be done.  

Since $C$ is a smooth curve at $P$, the ideal $P=\left<z_0,\ldots,z_{m-1},z_m-\gamma_j\right>$ must be principal.  On the curve, however, $z_0^{an}=z_m^a-1=\left(z_m-\gamma_j\right)\cdot p(z_m)$ where $p(\gamma_j)\neq0$, so that $z_m-\gamma_j$ fails to generate $P$.  Since $P$ is principal, this means that for some $0\leq i\leq m-1$ we have $P=\left<z_i\right>$ (i.e. $v_P\left(z_i\right)=1$).  Now consider the $i^{th}$ copy of equation \ref{eq:U1-0,1}, $z_1^i=z_iz_0^{i-\left\lceil\frac{ri}{m}\right\rceil}$ and take valuations of both sides.  Using equation \ref{eq:valrel-z0z1} to substitute for $v_P(z_1)$ one can solve for $v_P(z_i)$: \[ v_P(z_i)=\left(\left\lceil\frac{ri}{m}\right\rceil-\frac{ri}{m}\right)v_P(z_0).\]  Since the term in parentheses is a rational number whose denominator is a factor of $m$, the only way for $v_P(z_i)$ to be $1$ is for $v_P(z_0)$ to divide $m$.  

Consequently, $v_P(z_0)=m$, and $R_0=\sum_{j-1}^a(m-1)\cdot P$ concluding the proof.  
\end{proof}

Now that we have the degree and ramification divisor of $\phi:C\rightarrow{\mathbb P}^1$ in hand, determining $C$'s genus is a simple matter of plugging the results into Hurwitz's theorem.  

\begin{thm}\label{thm:WFCgenus}
The genus of a Fermat-type curve on ${\mathbb P}(1,m,n)$ of degree $amn$ is \[ g\left(C\right)=\frac{\displaystyle \left(am-1\right)\left(an-2\right)+a\left(m-1\right)}{\displaystyle 2}.\]
\end{thm}
\begin{proof}
Hurwitz's theorem \cite[Cor. IV.2.4]{Hart} states that given a finite map of curves, $f:X\rightarrow Y$ with ramification divisor $R$, the genus of each curve is related by \begin{equation}\label{eq:hurwitz} 2g(X)-2=\deg(f)\left(2g(Y)-2\right)+\deg(R).\end{equation}  Apply this result to the map $\phi:C\rightarrow{\mathbb P}^1$.  Lemma \ref{lem:degphi} gives $\deg(\phi)=am$.  In order to compute the degree of the ramification divsor, appeal to Lemmae \ref{lem:U0ramdiv} and \ref{lem:x0ramdiv} and the fact that $R=\bar{R}+R_0$ to see that \[ \deg(R)=\left(\sum_{j=1}^{an}am-1\right)+\left(\sum_{j=1}^am-1\right)=an(am-1)+a(m-1).\]  Since $g\left({\mathbb P}^1\right)=0$, equation \ref{eq:hurwitz} leaves \[ 2g(C)-2=am\left(-2\right)+an\left(am-1\right)+a\left(m-1\right) \] which may be solved for $g(C)$.  
\end{proof}


\begin{thebibliography}{amsplain}

\bibitem[D]{Dol} {\scshape Igor Dolgachev}, {\em Weighted Projective Varieties}. in {\em Group Actions and Vector Fields} (ed. J.B. Carrell), {\em Springer Lecture Notes in Math}. {\bf 956} (1982), 34--71.  \href{http://www.ams.org/mathscinet-getitem?mr=0704986}{MR0704986} (85g:14060). 
\bibitem[F]{Fult} {\scshape William Fulton}, {\em Introduction to Toric Varieties}.  Annals of Mathematics Studies, {\bf 131}.  The William H. Roever Lectures in Geometry. {\em Princeton University Press, Princeton, NJ}, (1993).  \href{http://www.ams.org/mathscinet-getitem?mr=1234037}{MR1234037} (94g:14028).  
\bibitem[H]{Hart} {\scshape Robin Hartshorne},  {\em Algebraic Geometry}.  Graduate Texts in Mathematics {\boldmath $52$}.  {\em Springer-Verlag, New York-Heideilberg}, (1977). \href{http://www.ams.org/mathscinet-getitem?mr=0463157}{MR0463157} (57 \#3116).  
\bibitem[K]{Ker} {\scshape Jeremiah M. Kermes}, {\em Desingularizations of Some Weighted Projective Planes}, submitted to Illinois J. Math., \href{http://arxiv.org/abs/0710.3409v1}{arXiv:0710.3409v1 [math.AG]} (2007).  
\bibitem[O]{Oda} {\scshape Tadao Oda}, {\em Convex Bodies and Algebraic Geometry -- An introction to the theory of toric varieties}. Translated from the Japanese. Results in Mathematics and Related areas (3), 15. {\em Springer-Verlag, Berlin}, (1988). \href{http://www.ams.org/mathscinet-getitem?mr=0922894}{MR0922894} (88m:14038).  

\end{thebibliography}
\end{document}